\documentclass[12pt,Preprint]{article}

\pdfoutput=1

\usepackage [english]{babel}
\usepackage[utf8]{inputenc}
\usepackage[autostyle]{csquotes}
\usepackage[T1]{fontenc}
\usepackage{amssymb}
\usepackage{amsmath}
\usepackage{amsfonts}
\usepackage{amsthm}
\usepackage{enumitem}
\usepackage[heavycircles]{stmaryrd}
\usepackage{caption}
\usepackage{subcaption}
\usepackage{mathrsfs}
\usepackage[colorlinks=true, linkcolor=blue]{hyperref}
\usepackage{cleveref}
\usepackage{verbatim}
\usepackage{graphicx} 
\usepackage[all,cmtip]{xy}

\newtheorem*{Theorem*}{Theorem}
\newtheorem{Theorem}{Theorem}[section]

\theoremstyle{definition}
\newtheorem{Definition}[Theorem]{Definition}

\newcommand{\aA}{\mathscr{A}}
\newcommand{\aB}{\mathscr{B}}

\newcommand{\aF}{\mathscr{F}}
\newcommand{\aG}{\mathbb{G}}
\newcommand{\asG}{\mathscr{G}}

\newcommand{\aO}{\mathscr{O}}

\newcommand{\am}{\mathfrak{m}}

\DeclareMathOperator{\aBr}{\operatorname{Br}}

\DeclareMathOperator{\aPic}{\operatorname{Pic}}

\DeclareMathOperator{\aSpec}{\operatorname{Spec}}

\DeclareMathOperator{\aker}{\operatorname{ker}}

\DeclareMathOperator{\fEnd}{\mathcal{E}\it{nd}}

\setlist[enumerate]{label=\bfseries(\roman*)}
\numberwithin{equation}{section}

\title{Pinching Azumaya algebras}

\author{Johannes Fischer}

\begin{document}

\maketitle

\begin{abstract}
We show, that for a morphism of schemes from X to Y, that is a finite modification in finitely many closed points, a cohomological Brauer class on Y is represented by an Azumaya algebra if its pullback to X is represented by an Azumaya algebra. Part of the proof uses an extension of a result by Ferrand, 
on pinching of finite locally free sheaves, to Azumaya algebras. 
\begin{footnotesize}  

\end{footnotesize}
\end{abstract}

\tableofcontents

\section*{Introduction}
\addcontentsline{toc}{section}{Introduction}

The classical theory of central simple algebras and Brauer groups over fields was generalized to schemes by Grothendieck \cite{Grothendieck:Le_groupe_de_Brauer_I-III}.
One important question of Grothendieck was whether the injective Brauer map $\delta: \aBr(X) \to \aBr'(X)$ from the Brauer group to the cohomological Brauer group is also surjective, i.e. if every cohomological Brauer class is represented by an Azumaya algebra.

In the general case the question is open to this day.
Grothendieck himself gave positive answers for schemes of dimension $\leq 1$ and for regular surfaces \cite{Grothendieck:Le_groupe_de_Brauer_I-III}. The perhaps best known result is that $\aBr(X)=\aBr'(X)$ holds for schemes that admit an ample invertible sheaf. The proof is due to Gabber, but was never published; there exists a proof by de Jong \cite{DeJong:A_result_of_Gabber}
(see also \cite{Colliot-Thelene_Skorobogatov:The_Brauer_Grothendieck_group} Theorem 4.2.1).
Further positive results are due to DeMeyer and Ford, for a smooth toric variety over an algebraically closed field of characteristic 0 \cite{DeMeyer_Ford:On_the_Brauer_group_of_toric_varieties}; and
Schröer, for a separated geometrically normal algebraic surface \cite{Schroer:There_are_enough_Azumaya_algebras_on_surfaces}. A recent paper of Mathur shows 
$\aBr(X)=\aBr'(X)$ for a separated surface \cite{Mathur:Experiments_on_the_Brauer_map_in_High_Codimension}.

A counterexample for a non-separated normal surface, where the Brauer map is not surjective,
was given by  Edidin, Hassett, Kresch and Vistoli (\cite{Edidin_Hasset_Kresch_Vistoli:Brauer_groups_and_quotient_stacks} Corollary 3.11).

The results in this paper are the main results of my PhD thesis \cite{Fischer:Pinching_Azumaya_algebras}. 
Let $f:X \to Y$ be a morphism of schemes. The pullback $f^*\aB$ of an Azumaya algebra $\aB$ on $Y$ is always an Azumaya algebra on $X$. I wondered what could be said in reverse. 
Let $\beta \in \aBr'(Y)$ such that $f^* \beta \in \aBr'(X)$ is represented by an Azumaya algebra $\aA$. Does there exist an Azumaya algebra $\aB$ on $Y$, with $f^*\aB \simeq \aA$? 
And does such a $\aB$ then also represent $\beta$? For the first question there are positive answers when $f$ is flat and surjective
(see for example \cite{Jahnel:Brauer_groups_Tamagawa_measures_and_rational_points_on_algebraic_varieties} Chapter III 2.8. Proposition). 

We can give a full answer in the following case:
\begin{Theorem*}(\Cref{Theorem: Pinch Azumaya algebra along finite modification in closed points})
Let $f:X \to Y$  be finite modification in finitely many closed points. A cohomological Brauer class $\beta \in \aBr'(Y)$
is represented by an Azumaya algebra on $Y$ if the pullback $f^*\beta$ is represented by an Azumaya algebra on $X$.
\end{Theorem*}
Such a morphism $f$ can for example arise by identifying two closed points in a scheme $X$, that have a common affine neighborhood, in order to obtain a scheme $Y$. The resolution of the singularity we just constructed is then a finite modification in one closed point $f:X \to Y$. By introducing mild singularities to a scheme $X$, where $\aBr(X)=\aBr'(X)$ is already known, the theorem can thus expand the class of schemes where we know that the Brauer map is surjective.

For example this gives a new scheme $Y$ with $\aBr(Y)=\aBr'(Y)$  if we identify two closed points on a smooth toric variety $X$ over a algebraically closed field of characteristic 0. Two closed points on such a toric variety do have a common affine neighborhood (\cite{Wlodarczyk:Embeddings_in_toric_varieties_and_prevarieties} Theorem A), so this identification indeed gives a scheme $Y$; but the resulting scheme is not normal, so not a toric variety and \cite{DeMeyer_Ford:On_the_Brauer_group_of_toric_varieties} does not apply any more. 
Another possible construction is to choose $X$ as the disjoint union of two schemes where $\aBr=\aBr'$ holds and glue them together in a closed point on each.

The proof works as following. In the setting of the theorem the conductor ideal $\mathcal{C}=\operatorname{Ann}_{\aO_Y}(f_* \aO_{X} / \aO_Y)$ determines closed embeddings both in $X$ and $Y$. Those in turn determine a commutative square that is both cartesian and cocartesian, called the conductor square. 
Given an Azumaya algebra $\aA$ on $X$ that represents $f^* \beta$, we extend a theorem of Ferrand \cite{Ferrand:Conducteur_descente_et_pincement} to Azumaya algebras, which allows us to 
construct an Azumaya algebra $\aB$ on $Y$, with $f^*\aB \simeq \aA$. We then use étale cohomology and the unique structure of the conductor square, to show that this $\aB$ represents $\beta$.

Given a cartesian and cocartesian square of schemes
\begin{equation*}
\begin{gathered}
\xymatrix{
 X' \ar[r]^-w \ar[d]_-g
& X \ar[d]^-f \\
Y' \ar[r]_-u
& Y,
}
\end{gathered}
\end{equation*}
where $f$ is affine and $u$ a closed immersion, Ferrand constructs an equivalence of categories for finite locally free sheaves
$$
	\Psi:  \operatorname{Loc}(Y') \times_{\operatorname{Loc}(X')} \operatorname{Loc}(X)
	 \longrightarrow
	\operatorname{Loc}(Y).
$$
Here the left hand side denotes the fiber product category constructed via pullbacks. 
We do the same construction for Azumaya algebras and show:
\begin{Theorem*}(\Cref{Theorem: Pinching Azumaya algebras})
The functor
$$
	\Psi:  \operatorname{Az}(Y') \times_{\operatorname{Az}(X')} \operatorname{Az}(X)
	 \longrightarrow
	\operatorname{Az}(Y)
$$
is an equivalence of categories. Furthermore, a quasicoherent $\aO_Y$-algebra $\aB$ is an Azumaya algebra on $Y$ if and only if
$u^* \aB \in  \operatorname{Az}(Y')$ and $f^*\aB \in  \operatorname{Az}(X)$.
\end{Theorem*}

This paper has two sections. The first section contains the proof of the main theorem. The second explains how to extend Ferrand's result to Azumaya algebras.

\

\textbf{Acknowledgments.} This article consists of parts of my PhD thesis. I wish to express my gratitude 
to my advisor Professor Dr. Stefan Schröer for giving me the opportunity to purse this project and for his continued advise while doing so. I would also like to thank the referee for his very useful comments, which helped to improve the paper.
This research was conducted in the framework of the research training group GRK 2240: Algebro-Geometric Methods in Algebra, Arithmetic and Topology, which is funded by the DFG.

\section{Pinching along finite modifications in closed points}
Let us shortly recall the definition for Azumaya algebras. Throughout, an algebra means a non zero, associative algebra with identity.
\begin{Definition}\label{Definition: Azumaya algebra}
Let $X$ be a scheme and $\aA$ and $\aO_X$-algebra, that is finite locally free as an $\aO_X$-module. We call $\aA$ an \textit{Azumaya algebra}, when the canonical homomorphism $\aA \otimes_{\aO_X} \aA^\circ \to \fEnd_{\aO_X\text{-mod}}(\aA)$ is an isomorphism.
\end{Definition}
Here $\aA^\circ$ denotes the opposite algebra. We will use in the next section that this condition is equivalent to 
$\aA_x \otimes_{\aO_X} \kappa(x)$ being a central simple algebra over the
residue field $\kappa(x)$ for all $x \in X$  (\cite{Colliot-Thelene_Skorobogatov:The_Brauer_Grothendieck_group} Theorem 3.1.1). Note that over fields central simple algebras are equivalent to Azumaya algebras.
We denote the Brauer map, from the Brauer group to the cohomological Brauer group,  by $\delta: \aBr(X) \to \aBr'(X)$. We say that an Azumaya algebra $\aA$ \textit{represents} a cohomological Brauer class $\alpha$, if $\delta([\aA])=\alpha$.

\begin{Definition}\label{Definition: finite modification}
We call a morphism of schemes  $f:X \to Y$ a \textit{finite modification in finitely many closed points} if the following holds true:
\begin{itemize}
\item	$f$ is finite.
\item  There exist dense open subsets $U \subset X$ and $V \subset Y$ such that $f(U) \subset V$ and $f_{\mid U}:U \to V$ is an isomorphism.
\item The image of $f$ is schematically dense.
\item $Y'=Y \setminus V$ shall consists of only finitely many closed points.
\end{itemize}
\end{Definition}
If the first three conditions hold, we call $f$ a \textit{finite modification}. Recall that the scheme theoretic image of $f$ is the closed subscheme defined by the ideal 
$\aker (\aO_Y \to f_* \aO_X)$. The scheme theoretic image is schematically dense in $Y$ if it is equal to $Y$; this is the case
if and only if $\aO_Y \to f_* \aO_X$ is injective.

Let $f:X \to Y$ be a finite modification.
We can define the \textit{conductor ideal}
$\mathcal{C}=\operatorname{Ann}_{\aO_Y}(f_* \aO_{X} / \aO_Y).$
The conductor ideal defines a closed subscheme $Y' \subset Y$. The inverse image $f^{-1}(\mathcal{C}) \aO_X$ defines
a closed subscheme $X' \subset X$, which is the base change of the closed immersion $u: Y' \to Y$; we have
$X'=Y' \times_Y X$.
So we get a cartesian square of schemes 
\begin{equation*}
\begin{gathered}
\xymatrix{
 X' \ar[r]^-w \ar[d]_-g
& X \ar[d]^-f \\
Y' \ar[r]_-u
& Y,
}
\end{gathered}
\end{equation*} 
It is also cocartesian. In the affine case this follows from \cite{Ferrand:Conducteur_descente_et_pincement} Lemme 1.2 and Théorème 5.1. Since every morphism in the square is affine this can be generalized with \cite{Ferrand:Conducteur_descente_et_pincement} Scolie 4.3.
Such a construction is called a \textit{conductor square}.

\begin{Theorem} \label{Theorem: Pinch Azumaya algebra along finite modification in closed points}
Let $f:X \to Y$  be finite modification in finitely many closed points. A cohomological Brauer class $\beta \in \aBr'(Y)$
is represented by an Azumaya algebra on $Y$ if the pullback $f^*\beta$ is represented by an Azumaya algebra on $X$.
\end{Theorem}

\begin{proof}
In the following all cohomology groups will be étale. We may assume that $Y$, and thus also $X$, have only finitely many connected components. 
This is possible by restricting to the connected components of $Y$ 
that contain one or more of the finitely many closed points in which we modify, since over any other connected component $f$ is an isomorphism.

Define a conductor square as above. Since $Y'$ is a scheme consisting of only finitely many closed points, we have $Y'=\aSpec(B)$, where $B$ is an Artin ring. And $X'=\aSpec(A)$ is the spectrum of an Artin ring $A$.
The resulting conductor square is:
\begin{equation}\label{Equation: square of schemes main theorem}
\begin{gathered}
\xymatrix{
 \aSpec(A) \ar[r]^-w \ar[d]_-g
& X \ar[d]^-f \\
\aSpec(B) \ar[r]_-u
& Y,
}
\end{gathered}
\end{equation}
with $h=ug=fw$.
Furthermore, we have a commutative square 
\begin{equation}\label{Equation: Brauer square main theorem}
\begin{gathered}
\xymatrixcolsep{4pc}
\xymatrix{
\aBr(B) \oplus \aBr(X) \ar[r]^-{\delta_B \oplus \delta_X}
& \aBr'(B) \oplus \aBr'(X) \\
\aBr(Y) \ar[r]_-{\delta_Y} \ar[u]^-{( u^*,f^*)}
& \aBr'(Y).  \ar[u]_-{( u^*,f^*)}
}
\end{gathered}
\end{equation}
The injective horizontal maps are given by the Brauer maps, while the vertical arrows come from the pullback maps for Brauer groups and cohomological Brauer groups, respectively.

For the first part of the proof we show that the pullback map 
$$( u^*,f^*):  \aBr'(Y) \to \aBr'(B) \oplus\aBr'(X)$$ is injective. This implies, that every morphism in the commutative square is injective.

Let $\varphi:T \to Y$ be a finite morphism of schemes.
The pullback map $\varphi^*:\aBr'(Y) \to \aBr'(T) $ factors as
$$\varphi^*:H^2(Y,\aG_{m}) \overset{\varphi_1}{\longrightarrow} H^2(Y,\varphi_*\aG_{m,T}) 
\overset{\varphi_2}{\longrightarrow}
H^2(T, \aG_{m,T})$$
(see \cite{Colliot-Thelene_Skorobogatov:The_Brauer_Grothendieck_group} 2.2.4(2)).
We show that
\begin{gather*}
( u^*,f^*):H^2(Y,\aG_{m}) \xrightarrow{(u_1,f_1)}
H^2(Y,u_*\aG_{m,B}) \oplus H^2(Y,f_*\aG_{m,X})
\\
\xrightarrow{(u_2,f_2)}
H^2(B, \aG_{m,B}) \oplus H^2(X, \aG_{m,X})
\end{gather*}
is injective, by showing that $(u_1,f_1)$ and  $(u_2,f_2)$  are injective.

The Leray-Serre spectral sequence (see \cite{stacks_project} Tag 03QC) gives the five term exact sequence
\begin{gather} \label{five term exact sequence}
\begin{gathered}
0 \longrightarrow H^1(Y, \varphi_* \aG_{m,T}) \longrightarrow H^1(T, \aG_{m,T})
\longrightarrow H^0(Y,R^1 \varphi_*  \aG_{m,T})
\\
\longrightarrow H^2(Y, \varphi_*  \aG_{m,T}) \overset{\varphi_2}{\longrightarrow} H^2(T, \aG_{m,T}).
\end{gathered}
\end{gather}
In this sequence is $ H^0(Y,R^1 \varphi_*  \aG_{m,T})=0$, since $R^1 \varphi_*  \aG_{m,T}=0$, for a finite morphism $\varphi$ (\cite{stacks_project} Tag 03QP). 
This implies that $\varphi_2$ is injective.
In particular, for $T=Y$ and $\varphi = u$ this shows that $u_2$ is injective, and for $T=X$ and $\varphi = f$ also that $f_2$ is injective. Thus $(u_2,f_2)$ is injective.

Since the conductor square (\ref{Equation: square of schemes main theorem}) is cartesian as well as cocartesian 
we have a short exact sequence:
\begin{equation*}
0 \longrightarrow  \aO_Y \longrightarrow  u_*\aO_{B} \oplus f_* \aO_{X} \longrightarrow h_* \aO_{A} \longrightarrow 0.
\end{equation*}
This exact sequence induces a short exact sequence
\begin{equation*}
1 \longrightarrow \aO_Y^\times  \longrightarrow  u_*\aO_{B}^\times \oplus f_* \aO_{X}^\times \longrightarrow h_* \aO_{A}^\times \longrightarrow 1.
\end{equation*}
Since base change by any étale morphism $U \to Y$ is flat, the exact sequence of Zariski sheaves above induces a short exact sequence of étale sheaves 
\begin{equation*}
1 \longrightarrow \aG_{m}  \longrightarrow  u_*\aG_{m,B} \oplus f_* \aG_{m,X} \longrightarrow h_* \aG_{m,A} \longrightarrow 1.
\end{equation*}
Taking cohomology we obtain an exact sequence
$$
 H^1(Y,h_* \aG_{m,A}) \longrightarrow
  H^2(Y,\aG_{m}) \longrightarrow  
  H^2(Y,u_*\aG_{m,B}) \oplus H^2(Y,f_*\aG_{m,X}).
$$
Using the five term exact sequence (\ref{five term exact sequence}) for 
$h: \aSpec(A) \to Y$ we see that 
$H^1(Y, h_* \aG_{m,A}) \simeq H^1(\aSpec(A), \aG_{m,A})$. Now 
$H^1(\aSpec(A), \aG_{m,A})=\aPic(A)$ and for an Artin ring $\aPic(A)=0$. This shows that 
$$(u_1,f_1):H^2(Y,\aG_{m}) \longrightarrow  H^2(Y,u_*\aG_{m,B}) \oplus H^2(Y,f_*\aG_{m,X})$$ 
is injective. Note that neither $u_1$ nor $f_1$ are necessarily injective; thus the same holds for $u^*$ and $f^*$

For the second part of the proof, let $\beta \in \aBr'(Y)$ be a cohomological Brauer class and  $\aA$ an Azumaya algebra on $X$ so that $\delta_X([\aA]) =f^* \beta$. Since 
$\aSpec(B)$ is affine there exist an Azumaya algebra $\aB'$ on it with $\delta_B([\aB']) =u^* \beta$.
Then  $(u^* \beta,f^* \beta) \in \aBr'(B) \oplus \aBr'(X)$ is in the image of $\delta_B \oplus \delta_X$
in the commutative diagram (\ref{Equation: Brauer square main theorem}). Since all maps in the diagram are injective it is enough to
construct an Azumaya algebra $\aB$ on $Y$ with $[u^* \aB] = [\aB']$ and $[f^* \aB] = [\aA]$. For then
$\delta_Y([\aB])=\beta$.

We will show in the next section, that if there exist an isomorphism $\tau: g^* \aB' \simeq w^* \aA$
we can pinch these Azumaya algebras together and obtain an Azumaya algebra $\aB$ on $Y$ with the required properties (\Cref{Theorem: Pinching Azumaya algebras}).

We still have to show that $\tau$ exists. By tensoring with appropriate matrix algebras we can replace $\aA$ and $\aB'$ by equivalent Azumaya algebras on $\aSpec(B)$ and $X$, that do have the same constant rank. 
Since $g^* \aB'$ and $w^* \aA$
both represent $h^* \beta$, they also have the same Brauer class in $\aBr(A)$.  
As a consequence of Wedderburn's theorem 
(see \cite{Gille_Szamuely:Central_simple_algebras_and_Galois_cohomology} Theorem 2.1.3) two Azumaya algebras over a field, that have the same class and the same rank, are isomorphic. 

It is left to show that $\aBr(A)$ is determined by Brauer classes of fields.
The structure theorem for Artin rings (\cite{Atiyah_Macdonald:Introduction_to_commutative_algebra}
Theorem 8.7) tells us that an Artin ring $A$ is uniquely, up to isomorphism, a finite direct product of local Artin rings $A_i$. It is $A=\prod_{i=1}^l A_i$. Now a local Artin ring is a henselian local ring and for such a rings holds $\aBr(R)=\aBr(R/ \am)$, where $\am$ denotes the maximal ideal (\cite{Milne:Etale_Cohomology} Chapter IV Corollary 2.13). Thus
$\aBr(A)=\bigoplus_{i=1}^l \aBr( A_i)=\bigoplus_{i=1}^l \aBr(k_i)$,
where the $k_i$ are the residue fields of the $A_i$.

\end{proof}
We remark that it should in principle be possible to expand this theorem to algebraic spaces.

\section{Pinching Azumaya algebras}

The following construction is due to Ferrand \cite{Ferrand:Conducteur_descente_et_pincement}. 
We start with a commutative square of schemes that is cartesian as well as cocartesian
\begin{equation}\label{Equation: commutative square of schemes}
\begin{gathered}
\xymatrix{
 X' \ar[r]^-w \ar[d]_-g
& X \ar[d]^-f \\
Y' \ar[r]_-u
& Y,
}
\end{gathered}
\end{equation} 
with $h=ug=fw$, where $f$ is affine and $u$ a closed immersion. Then $g$ is also affine and $w$ a closed immersion. In this case the 
pushout scheme $Y$ can be seen as the pinching of $X$ in the closed subscheme $X'$ along $g$. Also note that the conductor square from the previous section fulfills this conditions.
 
We denote the category of quasicoherent sheaves with $\operatorname{QCoh}(X)$. 
We can define a fiber product category  
$\operatorname{QCoh}(Y') \times_{\operatorname{QCoh}(X')} \operatorname{QCoh}(X)$ via pullbacks. The objects of the fiber product category are triples 
$(\asG', \tau, \aF)$, with quasicoherent sheaves $\asG' \in \operatorname{QCoh}(Y')$,
$ \aF \in \operatorname{QCoh}(X) $ and $\tau:g^*\asG' \to w^*  \aF$  is an isomorphism. The fiber product category comes with canonical forgetful functors, such that the diagram 
\begin{equation*}
\xymatrix{
\operatorname{QCoh}(X')  
&
 \operatorname{QCoh}(X) \ar[l]_-{w^*}
\\  
\operatorname{QCoh}(Y')  \ar[u]^-{g^*}
& 
\operatorname{QCoh}(Y') \times_{\operatorname{QCoh}(X')} \operatorname{QCoh}(X) \ar[u] \ar[l]
}
\end{equation*}
commutes up to a isomorphism of functors.

By the universal property of pullbacks of quasicoherent sheaves (basically the universal property of tensor products), the  equality $h=ug=fw$ in (\ref{Equation: commutative square of schemes}) induces a unique isomorphism of functors 
$$
\sigma: g^* u^*  \overset{\simeq}{\longrightarrow} w^*f^*.
$$
Note that $\sigma$ consists of a family of isomorphisms
$$
\sigma_{\asG}:g^* u^* \asG \overset{\simeq}{\longrightarrow} w^*f^* \asG
$$ 
for each quasicoherent $\aO_Y$-module $\asG$. One often thinks of this isomorphism of functors as "the identity", but to define the following functor into the fiber product category, we have to view it as a family of isomorphisms.
Then the universal property of fiber products, together with this data, gives a covariant functor
$$
\Phi:\operatorname{QCoh}(Y)  \longrightarrow     \operatorname{QCoh}(Y') \times_{\operatorname{QCoh}(X')} \operatorname{QCoh}(X), 
\quad
\asG \longmapsto (u^* \asG,\sigma_{\asG}, f^* \asG).
$$

We  define a right adjoint functor
$$
	\Psi:  \operatorname{QCoh}(Y') \times_{\operatorname{QCoh}(X')} \operatorname{QCoh}(X)
	 \longrightarrow
	\operatorname{QCoh}(Y).
$$
This functor $\Psi$ should map a triple  
$$(\asG', \tau, \aF) \in \operatorname{QCoh}(Y') \times_{\operatorname{QCoh}(X')} \operatorname{QCoh}(X)$$ to the fiber product 
$\Psi(\asG', \tau, \aF)=u_*\mathscr{\asG'} \times_{h_* \tau} f_* \aF$, so that
the diagram
\begin{equation}\label{Equation: sheaf functor cartesian}
\begin{gathered}
\xymatrix{
\Psi(\asG', \tau, \aF) \ar[rr] \ar[d] &
& 
f_*\aF \ar[d] 
\\
u_*\asG'\ar[r]
& 
u_* g_* g^*\asG' \ar[r]_-{h_*\tau} & f_* w_* w^*  \aF
}
\end{gathered}
\end{equation}
is cartesian. 
This diagram arises as follows. 
In general, for any morphism of schemes $\varphi: S \to T $ and any $\aO_T$-module $\aF$ the adjunction between pushforward and pullback gives a canonical morphism $\aF \to \varphi_* \varphi^* \aF$.
Since
$\tau:g^*\asG' \to w^* \aF$ is an isomorphism, the same holds for the pushforward
$h_*\tau: h_*g^*\asG' \to h_*w^* \aF$, which then allows us to define the fiber product.
We set the quasicoherent $\aO_Y$-module $\asG=\Psi(\asG', \tau, \aF)$ as the kernel of
$$
u_*\asG' \oplus f_*\aF 
\longrightarrow
h_* g^*\asG'  \oplus h_* w^*  \aF
\xrightarrow{h_* \tau - \operatorname{id}} h_* w^*  \aF.
$$
This gives a well defined fiber product. We pinched two quasicoherent sheaves together. Abusing notation we can now write 
$$\asG=u_*\asG' \times_{h_* \tau} f_* \aF.$$

Ferrand then restricts to the category of finite locally free sheaves and shows that in this case 
$\Psi$ 
is an equivalence of categories  (\cite{Ferrand:Conducteur_descente_et_pincement} Théorème 2.2 iv). The proof works by first showing that the functors $\Psi$ and $\Phi$ do map finite locally free shaves to finite locally free sheaves and then shows that $\Phi$ is the inverse of $\Psi$.
Note that Ferrand only treats the affine case. Since all morphisms in  (\ref{Equation: commutative square of schemes}) are affine the result generalizes.

We extend this result to Azumaya algebras. Denote the category of Azumaya algebras by $ \operatorname{Az}(X)$ and restrict the functors $\Psi$ and $\Phi$ to these categories.

\begin{Theorem}\label{Theorem: Pinching Azumaya algebras}
The functor
$$
	\Psi:  \operatorname{Az}(Y') \times_{\operatorname{Az}(X')} \operatorname{Az}(X)
	 \longrightarrow
	\operatorname{Az}(Y)
$$
is an equivalence of categories. Furthermore, a quasicoherent $\aO_Y$-algebra $\aB$ is an Azumaya algebra on $Y$ if and only if
$u^* \aB \in  \operatorname{Az}(Y')$ and $f^*\aB \in  \operatorname{Az}(X)$.
\end{Theorem}
\begin{proof}
We recall \Cref{Definition: Azumaya algebra} of Azumaya algebras. Let $\aB \in \operatorname{Az}(Y)$ and $(\aB', \tau, \aA) \in \operatorname{Az}(Y') \times_{\operatorname{Az}(X')} \operatorname{Az}(X)$. 
Since Azumaya algebras are finite locally free sheaves we know from Ferrand's result (\cite{Ferrand:Conducteur_descente_et_pincement} Théorème 2.2) that we have natural isomorphisms of functors, such that 
\begin{equation*} 
\Phi \Psi(\aB', \tau, \aA) \simeq (\aB', \tau, \aA) \quad \text{ and } \quad \Psi \Phi (\aB) \simeq \aB
\end{equation*}
are isomorphic as finite locally free sheaves. By construction the functors respect tensor products, so this 
morphisms respect algebra structures. If $\Psi$ and $\Phi$ map Azumaya algebras to Azumaya algebras, as claimed, they are inverse to each other.

The pullback of an Azumaya algebra is again an Azumaya algebra, so
$\Phi(\aB)=(u^* \aB,\sigma_{\aB}, f^* \aB)$ is an element of $\operatorname{Az}(Y') \times_{\operatorname{Az}(X')} \operatorname{Az}(X)$.
We will show that $\Psi(\aB', \tau, \aA) \in \operatorname{Az}(Y)$. From Ferrand's result we already know that it is a finite locally free sheaf.

The pushforward of a quasicoherent sheaf along an affine morphism is a quasicoherent sheaf. And
the quasicoherent sheaf $f_* \aA$ is an $\aO_Y$-algebra via $\aO_Y \to f_* \aO_X \to f_* \aA$.
So  $u_* \aB'$ and $f_* \aA$ are both $\aO_Y$-algebras, and $h_*\tau$ is an isomorphism of $\aO_Y$-algebras between them. 
Then we can view
$\Psi(\aB', \tau, \aA)=u_*\aB' \times_{h_* \tau} f_* \aA$
as a fiber product of $\aO_Y$-algebras, with the induced algebra structure.
This gives us a canonical way to equip the finite locally free sheaf $\Psi(\aB', \tau, \aA)$ with an $\aO_Y$-algebra structure. 

It is left is to show that the algebra $\bar{\aB}=\Psi(\aB', \tau, \aA)$ is actually an Azumaya algebra.
We have $\Phi(\bar{\aB}) \simeq (\aB', \tau, \aA)$. And thus $u^* \bar{\aB} \simeq \aB'$ and $f^* \bar{\aB} \simeq \aA$. Now $\aB'$ and $\aA$ are Azumaya algebras
and the algebra structure on $\bar{\aB}$ was constructed in such a way that this isomorphisms respect it.
This implies that $u^* \bar{\aB}$ and $f^* \bar{\aB}$ are Azumaya algebras. 
So for every 
$x \in X$,  $f^*(\bar{\aB}) \otimes_{\aO_{X}} \kappa(x) \simeq \bar{\aB}\otimes_{\aO_Y} \kappa(x) $ is an Azumaya algebra over $\kappa(x)$, with canonical isomorphism
$$
(\bar{\aB} \otimes_{\aO_{Y}} \bar{\aB}^{\circ}) \otimes_{\aO_Y} \kappa(x) \longrightarrow \fEnd_{\aO_Y \text{-mod}}(\bar{\aB}) \otimes_{\aO_Y} \kappa(x).
$$
Let $y=f(x)$, then the canonical homomorphism
$$
(\bar{\aB} \otimes_{\aO_{Y}} \bar{\aB}^{\circ}) \otimes_{\aO_Y} \kappa(y) \longrightarrow \fEnd_{\aO_Y \text{-mod}}(\bar{\aB}) \otimes_{\aO_Y} \kappa(y)
$$
is also an isomorphism, since $\kappa(x)/\kappa(y)$ is a field extension and thus faithfully flat. Thus $\bar{\aB} \otimes \kappa(y)$ is an Azumaya algebra. The same argument can be made for any $y \in Y$ that is in the image of $u$. Since the square is cartesian as well as cocartesian, every $y \in Y$ is in the image of either $u$ or $f$. This implies that $\bar{\aB} \otimes \kappa(y)$ is an Azumaya algebra for all $y \in Y$,
which shows that $\bar{\aB}$ is indeed an Azumaya algebra.
\end{proof}


\footnotesize

\textsc{Johannes Fischer, Mathematisches Institut, Heinrich-Heine-Universität, 40204 Düsseldorf, Germany}\par\nopagebreak
\textit{Email address:} \texttt{jfischer@math.uni-duesseldorf.de}

\end{document}